\def\today{\ifcase\month\or January\or February\or March\or April\or
    May\or June\or July\or August\or September\or October\or November
    \or December\fi\space\number\day, \number\year}
\magnification=1200
\baselineskip=8mm
\font\bigbf=cmb10 scaled \magstep2

\def\nz{\hfil\break\noindent}
\def\nl{\nz}
\def\lz{\vskip6mm\noindent}
\def\hz{\vskip7.2pt\noindent}
\def\shelahtitle{\centerline{\bigbf Large Normal Ideals Concentrating on
a
}
\centerline{\bigbf Fixed Small Cardinality}}
\font\msxm=msxm10 scaled 1200
\font\msym=msym10 scaled 1200
\def\vare{\varepsilon}
\def\name#1{\lower 1em\rlap{\char'176}#1}
\def\uhr{\hbox{\msxm\char'026}}
\def\eqdf{\buildrel\rm def\over =}
\def\lec{\mathrel{\le\mkern-11mu\raise1.85pt\hbox{$\scriptstyle\circ$}}}
\def\SI{{\cal I}}
\def\cf{{\rm cf}}
\def\otp{{\rm otp}}
\def\sup{{\rm  sup}}
\def\pr{{\rm pr}}
\def\levy{\,{\rm levy}\,}
\def\mod{{\rm mod}}

\def\rk{\rangle}
\def\lk{\langle}

\def\qed{\hfill\hbox{\msxm\char'003}}
\def\Vdash{\kern2pt{\vrule depth 1pt height .8 em width .4pt}\kern.1pt\vdash}
\def\force{\Vdash}
\def\notVdash{\;\>\rlap/\kern-5pt\Vdash}
\def\notvdash{\;\>\rlap/\kern-6pt\vdash}
\def\uhr{\hbox{\msxm\char'026}}
\def\name#1{\lower 1em\rlap{\char'176}#1}
\def\lq{{\rm``}}

\def\and{{\, \&\,}}

\def\smallbox#1{\leavevmode\thinspace\hbox{\vrule\vtop{\vbox
   {\hrule\kern1pt\hbox{\vphantom{\tt/}\thinspace{\tt#1}\thinspace}}
   \kern1pt\hrule}\vrule}\thinspace}
\font\msym=msym10
\textfont9=\msym

\def\levy{{\rm Levy}}
\def\eqdf{\buildrel\rm def\over =}
\def\SI{{\cal I}}
\font\funny=cmsy5
\def\supbarp{{\bar p}\llap{\lower3pt\hbox{$\textfont2=\funny\sim$\kern0.5pt}}{}}
\def\lec{\mathrel{\le\mkern-11mu\raise1.85pt\hbox{$\scriptstyle\circ$}}}
\today
\nl
 \shelahtitle
\bigskip

\centerline{\bf Saharon Shelah}
\medskip

\centerline{Institut of Mathematics}

\centerline{Hebrew University of Jerusalem}

\centerline{91904 Jerusalem, Israel}
\smallskip

\centerline{and}
\smallskip

\centerline{Department of Mathematics}

\centerline{Rutgers University}

\centerline{New Brunswick, NJ 08854, USA\footnote{}{Research supported
by ``Basic Research Foundation'' of The Israel Academy of
Sciences and Humanities.  Publication 542.}}
\bigskip
\bigskip
\bigskip

The property on the filter in Definition 1, a kind of large cardinal
property, suffices for the proof in
Liu Shelah [484] and is proved consistent as required there
(see conclusion 6).
A natural property which looks better, not only is not obtained
here, but is shown to be false (in Claim 7). On earlier related theorems
see Gitik Shelah [GiSh310].
\nl
\centerline{$*\qquad\qquad*\qquad\qquad*$}
\nl
{\bf 1. Definition}
(1) Let $\kappa$ be a cardinal and $D$ a filter
on $\kappa$ and $\theta$ be an ordinal $\le \kappa$ and $\mu<\chi$ but $\mu
\ge 2$ and $\chi\le \kappa$.
Let GM$_{\kappa, \chi, \theta, \mu}$ (D) be there following game:
\nl
a play lasts $\theta$ moves,
in the $\zeta's$ move the first player chooses a function $h_\zeta$ from
$\kappa$ to some ordinal $\gamma_\zeta<\chi$ and the second player
chooses a subset $B_\zeta$ of
$\gamma_\zeta$ of cardinality
 $<\mu$.

The second player wins a play if for every $\zeta<\theta$
the set
$\bigcap \{ \{\beta<\kappa: h_\varepsilon(\beta) \in B_\epsilon\}:
\epsilon\le \zeta\}$ is
$\not=\emptyset\,$mod$D$.
\nl
(2) If $\mu=2$ we may omit it,
if $\mu=2$ and $\chi=\kappa$ we omit $\chi$ and $\mu$.
\medskip\noindent
{\bf 2. Definition:}
 $(P\le, \le_{\pr})\in K_{\kappa, \chi,\theta, \mu}$ iff
\item{1.} $\kappa$ is a regular cardinal.
\item{2.} $(P,\le)$ is a forcing notion with minimal element
$\emptyset$ (if in doubt we use $\le_P$,
$\emptyset_P)$.
\item{3.} $P$ satisfies the $\kappa$-c.c.
\item{4.} $\le_{\pr}$ is a partial order on $P$ such that:
\itemitem{a.} $p\le_{\pr} q$ implies $p\le q$
\itemitem{b.} any $\le_{\pr}$-increasing chain of length $<\theta$
with first element $\emptyset$ in $P$ has an $\le_{\pr}$-upper bound.
\itemitem{c.} if $\gamma<\chi$ and $\name\tau$ is a $P$-name of an
ordinal $<\gamma$ and $\emptyset\le_{\pr} p\in P$ then for some $q$
and $B\subseteq \gamma$ of cardinality $<\mu$ we have
$p\le_{\pr} q\in P$ and $q$ forces $\name \tau\in B$.
\item{5.} for any $Y\subseteq P$ of cardinality $<\kappa$ there
is $P^*\lec P$ of cardinality $<\kappa$ such that
$P/P^*$ satisfies condition (4), i.e.
if $G^*\subseteq P^*$ is generic over $V$
and $P/G^*\eqdf\{p\in P:
p$ compatible with every $q\in G^*\}$
then
\itemitem{a.} in $P/G^*$,
any $\le_{\pr}$-increasing sequences starting with $\emptyset$ of length
$ <\theta$ have an $\le_{\pr}$-upper bound in $P/G^*$.
\itemitem{b.} if $p\in P/G^*$ and $\name\tau$ a $P$-name of an ordinal
$<\gamma$ where $\gamma<\chi$ then there is a subset $B$ of $\gamma$ of
cardinality $<\mu$ and $p'$,
$p\le_{\pr}$
$p'\in P/G^*$ such that $p'$
forces $\name\tau\in B$.
\hz
{\bf 2A Remark:} The relation in clause 4(b) is not really stronger than
having a winning strategy in the corresponding play, see [Sh250, 2.43]
(or [Sh-f, XIV 2.4]).
\medskip
\noindent
{\bf 3. Lemma:}
 Assume
\nl
a. $\kappa$ is a measurable cardinal with $D$ a
$\kappa$-complete ultrafilter on it
\nl
b. $(P\le, \le_{\pr})\in K_{\kappa, \chi, \theta, \mu}$
\nl
{\it Then}
in $V^P$ the second player wins GM$_{\kappa, \chi, \theta, \mu}$(D)
\nl
{\bf 3A Remark:}
1. We can replace ultrafilter by a filter in which the first player wins GM$_{\theta,
\kappa}$(D) [see Lemma 5].
\nl
{\it Proof:}
In $V$ we define a set $R$,
its members are sequences $\bar p=\langle p_\alpha:\alpha\in A^{\bar p}\rangle$
 where $A^{\bar p} \in D$ and $\emptyset\le_{\pr}p_\alpha
\in P$
(for $\alpha\in A^{\bar p}$).
On $R$ we define a partial order $\le_R$ as follows:
$\bar p \le_R \bar q$
iff
$A^{\bar q}\subseteq A^{\bar p}$
and for every $\alpha\in A^{\bar q}$ we have
$p_\alpha\le_{\pr} q_{\alpha}$.

Clearly, in $V$ the partial order $(R, \le_R)$ is $\theta$-complete.

For $G\subseteq P$ generic over $V$ we define $R[G]$ as
$\{\bar p:
\bar p\in R$ and
$\{\alpha\in A^{\bar p}:
p_\alpha\in G\}\not=\emptyset\,\mod D$
(in $V^P$,
$D$ is not a filter just a family of subsets
of $\kappa$ but it naturally generates a filter- just closed upward
and we refer to this filter in ``mod
$D$''$\}$.

For $G\subseteq P$ generic over $V$ and $\bar p\in R$ let
$w[\bar p, G]\eqdf\{\alpha \in A^{\bar p}:
p_\alpha\in G\}$.

So $R[G]=\{\bar p\in R:
w[\bar p, G]\not=\emptyset\,\mod D\}$. We now prove some facts.
\medskip
\noindent
{\bf 3B. Fact:}
In $V[G]$,
$(R[G], \le_R)$ is $\theta$-complete.
\nl
{\it Proof:}
If not then there is a $P$-name of a sequence of length $<\theta$,
 $\langle\name{\bar p}^\vare:\vare<\zeta\rangle$
and $r\in P$ which forces this sequence to be a counter example,
so $\zeta<\theta$.
So there are maximal antichains $\SI_\vare$ for $\vare<\zeta$
of conditions in $P$ forcing a value to $\name{\bar p}^\vare$ (note
$\name{\bar p}^\varepsilon$ is a $P$-name of a member of $V$);
let $Y$ be the set of elements appearing in some
$\SI_\varepsilon$ and $r$.
As $P$ satisfies the $\kappa$-c.c. clearly $Y$ has cardinality
$<\kappa$ so there is $P^*$ as required in condition (5) of Definition 2.
 Let $G^*\subseteq P^*$ be generic over $V$ and $r\in G^*$.

Now working in $V[G^*]$ we can (for each $\vare<\zeta$)
compute $\bar {\name p}^\vare $ and $A^{{\supbarp}^\vare}$,
 call it then $\bar p^\vare$ and
$A_\vare$ respectively and so $\bigwedge_\vare A_\vare\in D$ and
$A^*\eqdf \cap \{A_\vare:\vare<\zeta\}$
 belongs to $D^{V[G^*]}$
(=the ultrafilter which $D$ generates in $V[G^*]$, remember $\vert
P^*\vert<\kappa$, $D$ a $\kappa$-complete ultrafilter);
also letting $w_\vare\eqdf\{\alpha\in A^*:$
there is $G\subseteq P$ generic over $P$ extending $G^*$ to which
$p^\vare_\alpha$
belongs$\}\in V[G^*]$
we know that in $V[G]$ we get a $D$-positive set
$w[\name{\bar p}^\vare, G]$
 (because $r$ forces this) hence in
$V[G^*]$ the set $w_\vare$ is $D$-positive but in $V[G^*]$
we know $D^{V[G^*]}$ is an ultrafilter so necessarily $w_\vare$
 belongs to $D^{V[G^*]}$;
clearly for $\vare<\zeta$,
$\alpha\in w_\vare$ we have
$p^\vare_\alpha\in P/G^*$.
Let $B^*=A^*\cap \bigcap\{w_\vare:\vare<\zeta\}$,
it is in $D^{V[G^*]}$.
Now for any $\alpha \in B^*$ the sequence
$\langle p^\vare_\alpha:\vare<\zeta\rangle$
is a $\le_{\pr}$-increasing sequence of member of $P/G^*$
and by demand (5) (a) of Definition 2, the sequence
 has an
$\le_\pr$-upper bound $q_\alpha$ (in $P/ G^*$).
Let $r_\alpha \in G^*$ be above $r$ and force that
this holds and moreover force some specific $q_\alpha\in P_\alpha$ is as
above. So, still in $V[G^*]$, for some $C\in D$,
$C\subseteq B$ and
$r^*\in G^*$ we have
$(\forall\alpha\in C)[r_\alpha=r^*]$ without loss of generality $C\in V$.
As for $\alpha\in C\subseteq B$,
$r^*=r_\alpha\Vdash \lq q_\alpha\in P/\name G_{p^*}$'',
$r^*$ is compatible with every $q_\alpha$ $(\alpha\in C$).
By 3D below for some $q^+$,
$r^*\le q^+\in P$ and
 $q^+\Vdash_P\lq\{\alpha:q^+_\alpha\in \name G_{P^*}\}$''$
\neq \emptyset  \,\mod\, D$.
So $q^+$ (which is above $r\le r^*$) force that
 $\bar q=\langle q_\alpha: \alpha\in C\rangle$
is an upper bound as required.
(note: $\bar q\in V$, $r^*$ force it is an upper bound of
$\{\bar{\name p}^\vare:\vare<\zeta\}$;
we need $q^+_{\alpha(*)}$ as we do not know the value of
$\bar{\name p}^\vare$.
\qed$_{3B}$
\medskip
\noindent
{\bf 3C Fact:} Let $G\subseteq P$ be generic over $V$.
In $V[G]$,
if $\gamma<\chi$ and $\bar p \in R[G]$ and $h$ a function from
$\kappa$ to $\gamma$,
{\it then} for some $\bar q$ we have:

a. $\bar q \in R[G]$

b. $\bar p\le_R\bar q$

c. on $w[\bar p, G]$ the range of the function $h$ is of cardinality
$<\chi$.
\hz
{\it Proof:} Assume the conclusion fails then some $r\in G$ forces that
it fails for a
specific $\bar p$ and $P$-name
$\name h$
( so in particular $r$ forces that $w[\bar p, \name G]\not=\emptyset\, \mod D$.)
Let $w^*=:\{\alpha\in A^{\bar p}$: the conditions
$r, $ $p_\alpha$ are compatible in $P$
(equivalently,
$r$ does not force $\alpha\notin w[\bar p, \name G])\}$
(so $w^*\in V$) and $w^*\in D$.
 Now let $P^*$ be as in condition (5) of Definition 2
for $Y=\{r\}$
(so in particular $r\in P^*$).
Now:
\item{$(*)$} for every $\alpha\in w^*$ there are $r^*_\alpha$ and
$q_\alpha$ and
$B_\alpha$ such
that:
\itemitem{a.} $r\le r^*_\alpha\in P^*$.
\itemitem{b.} $p_\alpha\le_{\pr} q_\alpha$.
\itemitem{c.} $r_\alpha^*\force_{P^*} ``q_\alpha\in P/\name{G}_{P^*}$''.
\itemitem{d.} $q_\alpha$ forces (for $P$) that $\name h(\alpha)\in B_\alpha$
and for some set $B\subseteq \gamma$ $(B\in V)$, we have
$\vert B_\alpha\vert<\mu$.
\nl
[Why? for every $\alpha$ in $w^*$ we can find $G\subseteq P$ generic
over $V$ to which $r$ and $p_\alpha$ belong
(as $\alpha\in w^*$); hence $p_\alpha\in P/(G\cap P^*)$ hence some
$r_\alpha^*\in
G\cap P^*$ force this
(for $P^*$) so without loss of generality
$r\le r^*_\alpha$
(as $G\cap P^*$ is directed).
Now apply condition (5) of Definition 2 to
 $G\cap P^*$,
$p_\alpha$ and $\name h(\alpha)$
and we get some $B\subseteq \gamma$.
$\vert B\vert<\mu$ and
$q_\alpha\in P/(G\cap P^*$)
such that
 $p_\alpha\le_{\pr} q_\alpha\in P/(G\cap P^*)$
and $q_\alpha$ forces  $\name h (\alpha)\in B$.
Now increasing again $r^*_\alpha$ we get
$(*)$].

\noindent
 So we can find for $\alpha\in w^*$,
$r_\alpha, q_\alpha$ and $B_\alpha$ as in
$(*)$,
 (all in $V$);
let $A^* \subseteq w^*$
be such that $A^*\in D$ and
$\langle B_\alpha:\alpha\in w^*\rangle$ is constant on $A^*$
and also $r_\alpha$ is constantly $r^*$ (note: $D$ is $\kappa$-complete
$w^*\in D$, and $\kappa$ is strongly
inaccessible hence $\vert \gamma\vert^{<\mu}<\kappa$ and $\vert P^*\vert
<\kappa$.
Now some $q^+$, satisfying $r^*\le q^+\in P$,
 forces that $\lk q_\alpha: \alpha\in A^*\rk$
is in $R[G]$ by fact 3D below and so clearly is as required in the
Fact 3C.
\qed$_{3C}$
\medskip \noindent
{\bf 3D. Observation}
Assume $ \bar p=\langle p_\alpha: \alpha\in A\rangle$
$\in R$ and $r\in P$ is compatible
(in $P$) with every $p_\alpha$ (for $\alpha\in A$).
{\it Then} some $r^*$,
$r\le r^*\in P$,
force that $\bar p\in R[\name G_P$].
\hz
{\it Proof:}
Let $\SI$ be a maximal antichain of $P$ above $r$
 such that for every $q\in\SI$ we have
either $q\force_P$ ``$w[\bar p, \name G_P]$
is a subset of $A_q$'' where $A_q\subseteq \kappa$ and $\kappa\setminus
A_q\in D$
 {\it or} $q\force_P$ `` $w[\bar p, \name G_P]\not=\emptyset\, \mod D$.

{\it So} $\SI$ has cardinality $<\kappa$ and if the conclusion fails
then always the first possibility holds; now
 we let
$B\eqdf \bigcap\{\kappa\setminus A_q: q\in \SI\}$, clearly it belongs to $D$.
Now there is $\alpha\in B\cap A$ (as $B\cap A\in D)$ and
there is $r^*\in P$ above $r$ and above $p_\alpha$
(exist by assumption);
now $r^*$ force that $\alpha\in w[\bar p, \name G_p]\subseteq
A_q\subseteq \kappa\setminus B$, contradiction.
\qed$_{3D}$
\lz
{\bf 3E.}{\it Continuation of the Proof of Lemma 3:}
immediate for the Facts 3B, 3C.
\qed$_3$

Now we shall redo it all in another version:
\lz
{\bf 4. Lemma:}
(From Gitik [Gi] \S3, relaying on \S1 there, in different terminology).
Assume $\chi<\kappa$,
$\theta<\kappa$ a regular cardinal,
$\kappa$ is a measurable cardinal of order $\theta+1$ (i.e. there is a
coherent sequence of ultrafilters on $\kappa$ of length $\theta+1$,
see [Gi, \S3 p.293],
with $D$ an ultrafilter on $\kappa$ appearing in the $\theta$'th place
in the appropriate sequence.
  \nl
{\it Then} for some forcing notion $P$ we have
 \itemitem{(a)} $P$ of cardinality $\kappa$, $\force_P$ ``$\kappa$ is
strongly inaccessible''. 
\itemitem{(b)} $\{\delta:\, \force_P \lq\cf(\delta)=\theta$''$\}\in D$
\itemitem{(c)} $P\in K_{\kappa, \chi, \theta, 2}$ (in particular
$P$ satisfies the $\kappa$-c.c., $\leq_{\pr}$ for $P$ is called
$\leq_E$ in [Gi] (called Easton) 
\itemitem{(d)} For some $\le_{{\rm pr}}$ Condition (4) of Definition 2
is satisfied by $P$ (for $\mu=2)$. Moreover,
given any $\chi^*<\kappa$ and $Y\subseteq P$ of cardinality $<\kappa$ we can
find
$P^*\lec P$
as in clause (5) of Definition 2 replacing $\theta$ and $\chi$ by
$\chi^*$.
\medskip

\noindent
{\bf5. Claim:}
Under the assumptions of lemma 4,
if $\theta+\chi\le\mu=\cf (\mu)<\kappa$ let
$Q=P*(\levy(\mu, <\kappa))^{V^P}$
defining $(p_1, q_1)\le_{\pr}(p_2, q_2)$ iff
$p_1\le_P p_2$ and $p_2\force_P$ ``$q_1\le q_2\in\levy(\mu,
<\kappa)^{V^P}$''
\nl
{\it Then} $Q\in K_{\kappa, \chi, \theta, 2}$ and in $V^Q$, $\kappa=\mu^+=2^\mu$.
\hz
{\it Proof:} Easy.

\noindent
{\bf 5A Remark:} Actually in the conclusion of Claim 5 we can weaken
$\theta+\chi\le \mu$
to $\theta^++\chi\le \mu^+$ hence in the conclusion
$\chi=\mu^+(=\kappa$) is o.k.
This applies also to conclusion 6.

\noindent
{\bf 5B Remark:}
Of course Claim 5 and Definition 2 are formulated so that we get
consistency results justifying the name of the paper. We formulate below
 (conclusion 6) the one used in Liu Shelah [LiSh484].

\lz
{\bf 6. Conclusion:}
Assume $0=n_0<n_1<n_2<\dots<n_\ell, n_\ell+1<n_{\ell+1}$,
 $\kappa_{\ell+1}$ is a measurable of order
$\theta_\ell+1$ and for simplicity GCH holds and stipulate
$\kappa_0=\aleph_0$ and $\theta_{\ell+1}<\kappa_{\ell+1}$ is regular for
$\ell<\omega$, moreover $\theta_\ell\le
\kappa_{\ell+1}^{+(n_{\ell+1}-n_\ell)}$.

{\it Then} there is a forcing notion $P$ of cardinality $\le
2^{\Sigma_{\ell<\omega}
\kappa_{\ell+1}}$ which preserves
$\cf(\theta_{\ell+1})=\theta_{\ell+1}$, makes $\kappa_{\ell+1}$ to
$\aleph_{n_{\ell+1}}$ and preserves $(\kappa_\ell)^{+i}$ if
$i<n_{\ell+1}$, preserves
G.C.H. and for $\ell<\omega$ in $V^P$ the second player wins GM$_{\aleph_{\ell+1},
\aleph_{n_{\ell+1}-1}, \theta_{\ell+1},2}
(D_{\ell+1})$ for some $D_{\ell+1}\in V$, a normal ultrafilter on
$\kappa_{\ell+1}$ of order
$\theta_\ell+1$.
\hz
{\it Proof:} We use iteration $\langle P_i, Q_i:i<\omega\rangle$
described as follows:
$Q_{\ell}$=the forcing notion from lemma 5 (for
$\kappa=\kappa_{\ell+1}$,
$\theta=\theta_{\ell+1}$,
$\mu=\kappa^{+(n_{\ell+1}-n_\ell)}_\ell$ and
$\chi_{\ell+1}=\kappa^{+(n_{\ell+1}-n_\ell-1)}_{\ell+1})$,
the limit is a full support for pure extensions of the $\emptyset$ and finite
support otherwise (for the Levy collapse all conditions are
pure extensions of $\emptyset$). The checking is standard.
\hfill\qed$_{6}$
\nz
{\bf Discussion:} We shall now
prove that for a natural strengthening of Definition 2, we cannot get
consistency results. Specifically we cannot, in the game in Definition
2, let player I just decrease the present $D$-positive set.
\qed$_6$
\lz
{\bf 7. Definition:}
(1) Let $\kappa$ be a cardinal and $D$ a filter on $\kappa$ and
$\theta$ be an ordinal $\le \kappa$.
Let GM$^*_\theta(D)$ be the following game:

a play lasts $\theta$ moves;
 in the $\zeta$'s move 

{\it first player} chooses a subset $A_\zeta$
of $\kappa$,
$A_\zeta\not=\emptyset\,\mod D$ such that: if $\zeta=0$,
$A_\zeta\subseteq \kappa$ and if $\zeta=\varepsilon+1$ then
$A_\zeta\subseteq B_\varepsilon$ and if $\zeta$ is a limit ordinal
then $A_\zeta=\cap\{A_\vare:{\vare<\zeta}\}$

and then {\it the second player}  chooses a subset $B_\zeta$ of $A_\zeta$
satisfying
$B_\zeta\not=\emptyset\,\mod D$.

A player wins the play if he has no legal move (can occur only to the
first player in a limit stage),
if the play lasts $\theta $ moves then the second player wins.
\lz
{\bf 8. Definition:}
Let $\lambda$ be regular countable,
$S\subseteq \lambda$;
we say that there is a ($\le\theta$)-square for $S$ if:
there is a set $S^+$,
and sequence $\langle C_\alpha:\alpha\in S^+\rangle$ 
such that:
\nz
a. $S\subseteq S^+\subseteq \lambda$
\nz
b. for $\beta\in C_\alpha$
(so $\alpha\in S^+)$ we have:
$\beta\in S^+$ and $C_\beta=\beta\cap C_\alpha$.
\nz
c. $\otp(C_\alpha)\le \theta$ for $\alpha\in S^+$.
\nz
d. if $\delta\in S$ is a limit ordinal then $\delta=\sup(C_\delta)$
\nz
e. $C_\alpha$ is a closed subset of $\alpha$.
\lz
{\bf 9. Claim;}
1) Assume $\lambda $ is regular $>\theta$, $D$ is a normal
filter on $\lambda^+$ to which
$\{\delta:\cf(\delta)=\theta\}$ belongs.
{\it Then} in the game GM$^*_{\omega+1}(D)$
(see Definition 8 below)
the second player does not have a winning strategy.
\nl
2) Assume $\lambda$ is regular larger than
$\vert \theta\vert^+,$ $\theta$ an ordinal,
$D$ is a normal
filter on $\lambda$ to which a set $S$ belongs,
and for $S$ there is a $(\le \theta)$-square (as defined in Definition
7 above)
(or just every  $S\subseteq \lambda$,
$ S\not=\emptyset\,\mod D$
has a subset $S'$ for which there is a $(\le \theta$)-square.
$S'\not= \emptyset\,\mod D$).

{\it Then} in the game
GM$^*_{\omega+1}(D)$
(see Definition 8 below), the second player does not have a winning
strategy.
\lz
{\it  Proof:}
 Part (1) follows form part (2) as the assumption of part (2)
follows by [Sh 365, 2.14]
(or [Sh 351, Th. 4.1]).
So we concentrate on proving part (2).

So let $\langle C_\alpha:\alpha\in S^+\rangle$
be as in Definition 8.
So without loss of generality $S^+\in D$.
We divide $\{\delta:\delta<\lambda$,
$\cf(\delta)=\aleph_0\}$
to $
\vert \theta\vert ^+$ stationary sets
$\langle T_i:i<\vert\theta\vert^+\rangle$.
As $D$ is a normal ideal on $\lambda$,
$\vert\theta\vert^+<\lambda$, clearly
  for each stationary subset $S'$ of $S$ which is $D$-positive there
are $S^*\subseteq S'$ which is $D$-positive and
ordinal $j^*<\vert
\theta\vert^+$ such that for every $\alpha\in S^*$ we have: 
$C_\alpha\cup\{\alpha\}$ is disjoint to $T_{j(*)}$.

Now suppose the second player has a  winning strategy
in GM$^*_{\omega+1}(D)$ which we call Sty.
We can choose by induction on $n<\omega$ a sequence $\langle A_\rho, B_\rho,
\beta_\rho:\rho\in {}^n\lambda\rangle$ 
 such that

1. for every $\rho\in {}^n\lambda$ the sequence  $\langle A_{\rho\uhr k},
B_{\rho\uhr k}:k\le n\rangle$
 is an initial segment of a play of the game in which the second player
uses his winning strategy Sty

2. for some $j<\vert \theta\vert ^+$,
for every $\alpha\in A_{\langle \rangle}$ we have
$C_\alpha\cup\{\alpha\}$ is disjoint to $T_j$.

3. $\beta_\rho\in S^+$ and for every $\rho \in {}^n\lambda$ and $\alpha\in
A_\rho$ we have $\beta_\rho\in C_\alpha$.

4. for $\rho\in {}^n\lambda$ we have: $\beta_\rho$ is larger than
sup range $(\rho)$.

There is no problem to carry the definition (for clause (3) remember D
is a normal filter on $\lambda$);
now let
$E\eqdf \{\delta<\lambda:$
for every $\rho \in {}^{\omega>}\delta$ we have
$\beta_\rho <\delta\}$;
clearly $E$ is a club
 of $\lambda$ hence there is an ordinal $\delta\in E\cap T_j$;
so choose an increasing $\omega$-sequence $\rho$ of ordinals
$<\delta$ with limit $\delta$;
look at
 $\langle A_{\rho\uhr k}, B_{\rho\uhr k}: k< \omega\rangle$ which is an initial
segment of a play of the game in which the second player uses his
winning strategy
  Sty.
 Let now $B=\cap\{B_{\rho\uhr k}:k<\omega\}$;
if $\sup (B)>\delta$ (which holds if $B\not= \emptyset \mod
D)$, $\alpha\in B\setminus (\delta+1)$ then for every $n$,
$\beta_{\rho\uhr n}\in C_\alpha$. Note: as $\rho\in {}^\omega\delta$,
and $\delta\in E$ clearly $\beta_{\rho\uhr n}<\delta$;
so $\delta\ge \bigcup_{n<\omega}\beta_{\rho\uhr n}$;
as $\beta_{\rho\uhr n}\ge \sup{\,\rm range}(\rho\uhr n)$ necessarily
$\delta\le\bigcup_{n<\omega}\beta_{\rho\uhr n}$ so equality
holds.
Hence also $\delta=(\cup_{n<\omega}$
$\beta_{\rho\uhr n}) \in C_\alpha$ (as
$\alpha>\delta=\bigcup_{n<\omega}\beta_{\rho\uhr n}$). So  $\delta\in C_\alpha$
but $\delta\in T_j$ whereas $\alpha\in B_{\langle \rangle}$,
contradiction.
 So $B$ is a subset of $\delta+1$,
contradicting to ``Sty is a winning strategy''.
\lz
{\bf 9A Remark:} This continues the argument that e.g. not for every
stationary $S\subseteq \{\delta<\aleph_3:\cf(\delta)=\aleph_0\}$, there
is a club E of $\aleph_3$ such that $\delta\in E\,\&\,$ if
$(\delta)=\aleph_2\Rightarrow S\cap \delta$ stationary in $\delta$
(find pairwise disjoint
$S_i\subseteq \{\delta<\aleph_3:\cf(\delta)=\aleph_0\}$, for
$i<\aleph_3$, if for $S_i$ we have $E_i$, choose $\delta\in
\bigcap_{i<\aleph_2}E_i$ of cofinality $\aleph_2$.
\bigskip
\centerline{{\bf References}}
\item{[Gi]} M. Gitik. Changing cofinalities and the non stationary
ideal, {\it Israel Journal of Mathematics} 56(1986) pp.280-314 .
\item{[GiSh310]} M.Gitik and S.Shelah, Cardinal preserving ideals,
{\it Journal of Symbolic Logic}, to appear.
\item{[LiSh 484]}  K.Liu and S.Shelah,
Confinalities of elementary substructures of structures on
$\aleph_\omega$, {\it Israel Journal of Mathematics}, to appear.
\item{[Sh 351]}  S.Shelah, Reflecting stationary sets and successors of 
singular cardinals,  
{\it Archive fur Math Logic} 31 (1991) pp. 25--34.
\item{[Sh 365]} S.Shelah, There are Jonsson algebras in many 
inaccessible cardinals, Ch III {\it Cardinal Arithmetic, OUP}. accepted.
\bye